\documentclass[11pt]{article}
\usepackage{latexsym}
\usepackage{amscd}
\usepackage{amssymb}
\usepackage{amsmath}
\usepackage{amsthm}
\usepackage{enumerate}
\usepackage{verbatim}
\def\co{\colon\thinspace}
\def\d{\delta}
\def\e{\epsilon}
\def\a{\alpha}

\def\g{\gamma}
\def\s{\sigma}
\newtheorem{thm}{Theorem}[section]
\newtheorem{cor}[thm]{Corollary}

\newtheorem{lem}[thm]{Lemma}
\newtheorem{prop}[thm]{Proposition}

\newtheorem{Example}[thm]{Example}

\newtheorem{remark}[thm]{Remark}
\newenvironment{rmk}{\begin{remark}\rm}{\end{remark}}
\newtheorem{Fact}[thm]{Fact}

\newtheorem{Nothing}[thm]{$\!\!\!$}
\newenvironment{nothing}{\begin{Nothing}\rm}{\end{Nothing}}
\begin{document}
\abovedisplayskip=6pt plus3pt minus3pt
\belowdisplayskip=6pt plus3pt minus3pt
\title{On Mostow rigidity for variable negative curvature}
\author{Igor Belegradek}
\date{}
\maketitle
\begin{abstract} 
We prove a finiteness theorem for the class of complete finite volume 
Riemannian manifolds with pinched negative sectional curvature, fixed
fundamental group, and of dimension $\ge 3$. 
One of the key ingredients is that
the fundamental group of such a manifold does not admit a small nontrivial
action on an $\mathbf R$-tree.
\end{abstract}
\section{Introduction}
\label{Introduction}
According to the Mostow rigidity theorem, 
the isometry type of a complete finite volume locally
symmetric negatively curved Riemannian $n$-manifold with $n\ge 3$
is uniquely determined by its fundamental group.
This is no longer true for finite volume manifolds of variable
sectional curvature. In fact, there even exist negatively curved manifolds
which are homeomorphic but not diffeomorphic to 
finite volume hyperbolic manifolds~\cite{FJO}. 
Yet it turns out that there
are essentially finitely many possibilities for the 
geometry and topology
of such manifolds provided the sectional curvature
is pinched between two negative fixed constants.

For reals $a\le b\le 0$, let $\mathcal M_{a, b,\pi, n}$ be the class
of complete finite volume Riemannian manifolds of dimension $n\ge 3$
with sectional curvatures in $[a, b]$ and 
fundamental groups isomorphic to $\pi$.
Note that $n$ can be read off the fundamental group $\pi$,
namely, $n=\mathrm{max}(\mathrm{cd}(\pi), \mathrm{cd}(N)+1)$ 
where ``cd'' is the cohomological dimension
and $N$ is a maximal nilpotent subgroup of $\pi$.
Also, $\mathrm{cd}(\pi)=n$ iff each manifold in 
$\mathcal M_{a, b,\pi, n}$ is closed.
Here is our main result.

\begin{thm}\label{intro: main thm} 
The class $\mathcal M_{a, b,\pi, n}$ falls into
finitely many diffeomorphism types. Furthermore,
for any sequence of manifolds in $\mathcal M_{a, b,\pi, n}$,
there exists a subsequence $(M_k,g_k)$,
a smooth manifold $M$, and 
diffeomorphisms $f_k:M\to M_k$ 
such that the pullback metrics $f_k^\ast g_k$ 
converge in $C^{1,\alpha}$ topology
uniformly on compact subsets to a complete finite volume
$C^{1,\alpha}$ Riemannian metric on $M$. 
\end{thm}

Although this result does not appear in the literature, 
it is not new except in dimension $4$ where only homeomorphism 
finiteness has been known.
However, the proof we present is very different from the existing
argument that runs as follows.

M.~Gromov and W.~Thurston~\cite{Gro5} used
straightening and bounded cohomology
to deduce that volume is bounded by the
simplicial volume:
$$vol(M)\le C(a/b,n)||M||<\infty$$
for any $M\in \mathcal M_{a, b,\pi, n}$. 
Since all the manifolds in $\mathcal M_{a, b,\pi, n}$ 
belong to the same proper homotopy type, they must
have equal simplicial volumes. 
Thus, volume  is uniformly bounded from above
on $\mathcal M_{a, b,\pi, n}$. 
Furthermore, a universal lower volume bound 
comes from the Heintze-Margulis theorem~\cite{BGS}.
For closed manifolds of dimension
$n\ge 4$ and sectional curvatures within $[-1,0)$, 
the diameter can be bounded in terms of volume~\cite{Gro6}.
Hence the conclusion of~\ref{intro: main thm} follows
from the Cheeger-Gromov compactness theorem. 
Similarly, the conclusion of the theorem~\ref{intro: main thm}
for finite volume manifolds of dimension $\ge 5$ can be deduced from
the work of K.~Fukaya~\cite{Fuk2}. The dimension restriction comes from
treating the ends by the weak h-cobordism theorem. In fact, Fukaya
proves a similar statement for $n=4$ where diffeomorphism finiteness
is replaced by homotopy finiteness.
(The topological $4$-dimensional weak $h$-cobordism theorem, 
unavailable at the time of~\cite{Fuk2}, can be 
also applied here because we deal with virtually nilpotent fundamental 
groups~\cite{Gui, FT}. 
Yet this only gives homeomorphism rather than diffeomorphism
finiteness.) 

By contrast, main ideas in our approach come from Kleinian groups
and geometric group theory. Essentially, given a degenerating 
sequence of manifolds $M_k$, 
one can use rescaling in the universal covers
to produce a nontrivial action of $\pi$ on an $\mathbf R$-tree
with virtually nilpotent arc stabilizers (cf.~\cite{Bes, Pau1, Pau2}).
Then results of Rips, Bestvina and Feighn~\cite{BF}
imply that $\pi$ splits over a virtually nilpotent subgroup.
We prove that this does not happen if $\mathcal M_{a, b,\pi, n}$
is nonempty. Then the methods of the Cheeger-Gromov compactness theorem
imply that $M_k$ subconverges in pointed 
$C^{1,\alpha}$ topology to a complete 
$C^{1,\alpha}$ Riemannian manifold $M$.
We then prove that $\pi_1(M)$ contains a subgroup
isomorphic to $\pi$ which implies that $M$ has finite volume
and, in fact, $\pi_1(M)\cong\pi$.
Now the convergence
$M_k\to M$ is analogous to strong convergence of Kleinian groups.
Studying this convergence yields the theorem~\ref{intro: main thm}.
Instead of using the h-cobordism theorem to deal with ends,
we find a direct geometric argument that works in all
dimensions. 
Similarly to~\cite{Gro5}, our methods provide a uniform
upper bound on the volume of manifolds in $\mathcal M_{a, b, \pi, n}$.

Note that the 
real Schwarz lemma of Besson-Courtois-Gallot~\cite{BCG2} 
gives yet another way to get a uniform
upper bound on volume of {\it compact} manifolds in 
$\mathcal M_{a, b, \pi, n}$, and hence another proof of
~\ref{intro: main thm} in this case.
It is still an open question whether their method extends to
finite volume manifolds.

Technically, our proof is much easier for closed manifolds; 
this case was previously treated in~\cite{Bel4}. 
One of the key facts 
needed for the theorem~\ref{intro: main thm} is 
the following lemma which is of independent interest.

\begin{lem} \label{intro: no split}
If the class
$\mathcal M_{a, b, \pi, n}$ is nonempty, then
$\pi$ does not split over a virtually nilpotent subgroup. Furthermore,
$\mathrm{Out}(\pi )$ is finite and $\pi$ is cohopfian.
\end{lem}

That $\pi$ does not split over a virtually nilpotent subgroup can
be also deduced (after some extra work) from a recent paper of 
B.~Bowditch~\cite{Bow3} who studied the structure of the 
splittings of relatively hyperbolic groups over
subgroups of peripheral groups. 
Unlike Bowditch's work that applies in a more general situation,
our argument is completely elementary. 

Recall that a group $\pi$ is called cohopfian if it
has no proper subgroups isomorphic to $\pi$. 
That $\mathrm{Out}(\pi )$ is finite and $\pi$ is cohopfian
is due to G.~Prasad~\cite{Pra} in locally symmetric case. 
Here we follow the idea of F.~Paulin~\cite{Pau2}, and
E.~Rips and Z.~Sela~\cite{RS} who proved these properties for
a large class of word-hyperbolic groups.

Topology of complete finite volume negatively curved manifolds 
seems to be encoded in the fundamental group.
It is a deep recent result of F.~T.~Farrell and L.~Jones~\cite{FJ3, FJ6} 
that any homotopy equivalence of two manifolds
in the class $\mathcal M_{n,a,b,\pi}$ with $n\ge 5$
is homotopic to a homeomorphism (which is a diffeomorphism
away from compact subsets). 
Then the smoothing theory~\cite{HM, KS} implies that 
there exist at most finitely many nondiffeomorphic manifolds in the class
$\mathcal M_{a,b,\pi,n}$ with $n\ge 5$.
Furthermore, if $n\ge 6$, there are 
homeomorphic negatively curved manifolds 
that are not diffeomorphic~\cite{FJ1, FJ2, FJ7, FJO}.
Such examples are still unknown if $n=4,5$. 
If $n=3$, one expects that any two manifolds
in the class $\mathcal M_{a,b,\pi, n}$ 
are diffeomorphic. This is known for Haken manifolds~\cite{Wald}
(note that any noncompact finite volume manifold is Haken).
For non-Haken manifolds this would follow from
(as yet unproved) Thurston's hyperbolization conjecture
and Mostow Rigidity.

Note that~\ref{intro: main thm} gives diffeomorphism finiteness 
in all dimensions. As we explained above this is mostly interesting if $n=4$;
in fact, the interior of any compact smooth $4$-manifold with nonempty
boundary has at least countably many smooth structures~\cite{BE}.

Theorem~\ref{intro: main thm}  
combined with results of Gao~\cite{Gao}
immediately implies that the class of finite volume Einstein manifolds
that lie in  $\mathcal M_{a, b, \pi, n}$ is compact in pointed
$C^\infty$ topology
(i.e. for any sequence of Einstein manifolds in $\mathcal M_{a, b, \pi, n}$,
there is a subsequence $M_k$, a manifold $M$, and diffeomorphisms $f_k\co M\to M_k$
such that $f_k^\ast g_k$ converge in $C^\infty$ topology uniformly on compact subsets
to a complete finite volume Einstein metric on $M$). 
Furthermore, since negatively curved
Einstein metrics on compact manifolds of dimension $>2$ are isolated in the
moduli space of the Einstein metrics~\cite[12.73]{Bes}, 
we conclude that, up to homothety, 
there are only finitely many compact Einstein manifolds
in $\mathcal M_{a, b, \pi, n}$. 
It is a tantalizing open problem to decide
whether any compact negatively curved Einstein manifold is locally symmetric.
A related question was recently resolved in dimension four:
any Einstein metric on a compact negatively curved 
locally symmetric $4$-manifold is
locally symmetric~\cite{BCG1, LeB}

One can also use~\ref{intro: main thm} to deduce several pinching results. 
For example, given a group $\pi$, there exists an $\epsilon=\epsilon(\pi) >0$
such that any finite volume manifold from $\mathcal M_{-1-\epsilon,-1,\pi , n}$
is diffeomorphic to a real hyperbolic manifold.
Note that $\epsilon$ has to depend on the topology of the manifolds 
(in our case on the fundamental group)
as examples~\cite{FJ1,FJ7, FJO, GT} show.
Similar results hold for almost quarter pinched K\"ahler 
and quaternionic-K\"ahler manifold manifolds. 

The structure of the paper is as follows. 
In the section $2$ we give some
background on finite volume negatively curved manifolds.
Sections $3$ and $5$ contain a proof of~\ref{intro: no split}.
Convergence of finite volume negatively curved manifolds
is discussed in Section $4$. 
Theorem~\ref{intro: main thm} is proved in Section $6$.
Applications to pinching are discussed in Section $7$. 

The author is grateful to Christoph B\"ohm, Brian Bowditch, 
Harish Seshadri, Peter Shalen,  Christopher Stark, and
McKenzie Wang for helpful discussions or communications.

\section{Preliminaries}

Let $M$ be a finite volume complete Riemannian manifold
of sectional curvature within $[a,b]$, $a\le b<0$ of dimension $n>2$. 
In this section we list some properties of $M$ which 
we are going to use throughout the paper without explicit references. 
A comprehensive account on nonpositive curvature
can be found in~\cite{BGS}. 
Also see~\cite{HI, Sch, Bow2, Bow1}.

\begin{nothing}
\bf{Virtually nilpotent subgroups.}\rm\
Write $M$ as $X/\pi$ where $X$ is the universal cover of $M$
and $\pi\cong\pi_1(M)$ is the covering group.
By the Cartan-Hadamard theorem $X$ is diffeomorphic
to the Euclidean space, thus $M$ is aspherical. 
An infinite order isometry $\gamma$ of $X$ either stabilizes
a bi-infinite geodesic or fix exactly one point at infinity;
such a $\gamma$ is called {\it hyperbolic} or {\it parabolic},
respectively.

Virtually nilpotent discrete subgroups 
of $Isom(X)$ are finitely generated~\cite{Bow2}. 
Since $X$ is contractible, $\pi$ is torsion free.
Any nontrivial virtually nilpotent subgroup of $\pi$ is 
either an infinite cyclic group generated by a loxodromic isometry
and stabilizing a bi-infinite geodesic, or 
a group that consists of parabolic isometries with a common fixed point
at infinity. Any virtually nilpotent subgroup of $\pi$ lies in a
{\it unique} maximal virtually nilpotent subgroup of $\pi$.
See~\cite{BGS, Bow2, Bow1} for more information.
\end{nothing}

\begin{nothing}
\bf{Thin/thick decomposition.}\rm\
For a positive $\epsilon$, 
write $M_{[\e,\infty)}$ for the {\it $\e$-thick part} of $M$
which is the set of points of
$M$ with injectivity radius $\ge\epsilon$.
Similarly, $M_{(0,\e)}=M\setminus M_{[\e,\infty)}$ is called 
the {\it $\e$-thin part} of $M$.

According to the Margulis lemma, there exists a universal constant
$\mu_{n,a}$ such that for each $\epsilon<\mu_{n,a}$
the $\e$-thin part of $M$ is a union of {\it finitely many} connected components.
Unbounded components are called {\it cusps} while 
bounded components are called {\it tubes}. Each tube contains a closed
geodesic of length $\le 2\epsilon$ and is homeomorphic to the tubular
neighborhood of the geodesic. 
Since there are only finitely many tubes, we can assume that
$M_{(0,\epsilon)}$ consists of cusps by taking $\epsilon$ small enough.

Each cusp is a union of geodesic rays 
emanating from a common point at infinity.
Also, given a cusp $C$, let  $\tilde{C}$ be
a connected component of the preimage of $C$ in the universal cover. 
The group $\{\gamma\in\pi:\gamma(\tilde{C})=\tilde{C}\}$
coincides with the stabilizer $\Gamma_z$ in $\Gamma$
of a point $z$ at infinity. The group $\Gamma_z$ preserves
horospheres centered at $z$ and acts on each horosphere with
compact quotient since $C=\tilde{C}/\Gamma_z$ has finite volume.
Horospheres are $C^2$ submanifolds of $X$ each diffeomorphic to 
the Euclidean space~\cite{HI}. 
Each horosphere centered at $z$ is orthogonal to
geodesics asymptotic to $z$. Tangent vectors to
such geodesics form the so-called {\it radial} 
vector field on $X$ and $C$; this is a $C^1$-vector field~\cite{HI}.
(Throughout the paper all geodesic are assumed to have unit speed.)
In fact, radial vector field is the gradient of the so called
Busemann function. Any Busemann function defines a $C^2$-Riemannian
submersion of a cusp region bounded by a horosphere into the real line. 
Each cusp is diffeomorphic to
the product of a real line and a closed aspherical
manifold which is the quotient of a horosphere by $\Gamma_z$.
Note that the boundary $\partial C$ of a cusp is generally nonsmooth.
Pushing along geodesic rays asymptotic to $x$
defines a homeomorphism of $\partial C$ and the 
$\Gamma_z$-quotient of a horosphere.
\end{nothing}

\begin{nothing}
\bf{Compactification.}\rm\
According to~\cite{BGS}, $M$ is diffeomorphic to
the interior of a compact manifold 
with (possibly empty) boundary.
If the boundary is nonempty, its connected
components are quotients of horospheres by maximal
parabolic subgroups of $\pi$. 
Each boundary component corresponds 
to a conjugacy class of maximal parabolic subgroups of $\pi$.
The inclusion of each boundary component is $\pi_1$-injective.
\end{nothing}

\begin{nothing} \label{prelim: exp conv}
\bf{Exponential convergence of geodesics.}\rm\
To simplify notations we let $a=-K^2$, $b=-k^2$.
Let $\g_1$, $\g_2$ be two geodesics asymptotic to a point
$z=\g_1(+\infty)=\g_2(+\infty)$ such that $\g_1(0)$, $\g_2(0)$
lie on the same horosphere centered at $z$.
Then for any $t$, $\g_1(t)$, $\g_2(t)$ lie on the same horosphere 
centered at $z$. Denote by $h(t)$
the distance between $\g_1(t)$ and $\g_2(t)$
on this horosphere equipped with the induced Riemannian metric.
Also denote by $d(t)$ the distance
between $\g_1(t)$ and $\g_2(t)$ in $X$.
It is proved in~\cite{HI} that for $t\ge 0$
\[e^{kt}\le \frac{h(0)}{h(t)}\le e^{Kt}\ \ \ \mathrm{and}
\ \ \ \frac{2}{k}\sinh\left(\frac{kd(t)}{2}\right)\le h(t)\le 
\frac{2}{K}\sinh\left(\frac{Kd(t)}{2}\right).\]
Therefore, we deduce that $d(t)\le h(t)\le h(0)e^{-kt}\le  
\frac{2}{K}\sinh\left(\frac{Kd(0)}{2}\right)e^{-kt}$ and 
$d(0)e^{-Kt}\le h(0)e^{-Kt}\le h(t)\le 
\frac{2}{K}\sinh\left(\frac{Kd(t)}{2}\right)$ for all $t\ge 0$.
It is straightforward to check that $\sinh(x)\le 2x$ whenever
$x\in [0,1]$. Also, it follows from comparison with Euclidean triangles that
$d(t)\le d(0)$ for $t\ge 0$. Hence, if $d(0)\le 2/K$ and $t\ge 0$,
then $e^{kt}/2\le d(0)/d(t)\le 2e^{Kt}$.
\end{nothing}

\section{No splitting}

Throughout this section $\pi$ is the fundamental group 
of a finite volume noncompact complete Riemannian manifold $M$ of dimension
$n>2$ and with sectional curvatures within $[a,b]$ for $a\le b<0$. 
By~\cite{BGS} $M$ is the interior of a compact manifold
with boundary which we denote by $M_\pi$.

A group $G$ is said to {\it split over a subgroup} $C$ 
if $G=A\ast_C B$ or
$G=A\ast_C$ where $A\neq C\neq B$.
It is well-known that $A$ and $B$ necessarily have infinite index in $G$.
Note that $G$ splits over $C$ iff 
$G$ act without edge inversions on a simplicial
tree with no proper invariant subtree, no global fixed point, and exactly
one orbit of edges such that $C$ is a stabilizer of some edge~\cite{Ser}.

The purpose of this section is to prove that $\pi$
does not split over a virtually nilpotent group.
When $M$ is a closed manifold, 
this can be shown simply by looking at the Mayer-Vietoris
sequence of the splitting (see e.g.~\cite{Bel4}). 

Noncompact case is more subtle.
Main ingredients of the proof are again Mayer-Vietoris
sequence, and the following splitting-theoretic lemma of B.~Bowditch.
More general (and harder to prove) results
can be deduced from a recent paper of Bowditch~\cite{Bow3} 
who studies the structure of splittings of relatively 
hyperbolic groups over subgroups of peripheral groups. 
By contrast, our approach is elementary: we use basic manifold
topology and Bass-Serre theory.

\begin{lem} \label{relative split}
Suppose that $\pi$ splits over a virtually nilpotent subgroup.
Then $\pi$ also splits over a virtually nilpotent 
subgroup $C$ as $A\ast_C B$ or
$A\ast_C$ where $A\neq C\neq B$
in such a way that a conjugate of 
any maximal parabolic subgroup lies in $A$ or $B$.
Furthermore, if $C$ is parabolic, then the splitting can be 
chosen so that
the maximal parabolic subgroup containing $C$
lies in $A$.
\end{lem}
\begin{proof} 
The proof is essentially borrowed from~\cite[3.5, 5.2]{Bow3}
where a more general situation is considered.
We specialize terminology to our case and give more details
when seems appropriate.

Since $\pi$ splits over a virtually nilpotent subgroup,
$\pi$ acts without edge inversions on a simplicial
tree $T$ with no proper invariant subtree, no global fixed point, and exactly
one orbit of edges. 
We seek to construct a $\pi$-action on a (perhaps another) tree with
the above properties such that every maximal parabolic subgroup of $\pi$
fixes a vertex.

Fix an edge $e$ of $T$ and denote its stabilizer in $\pi$ by $\pi_e$.
Let $P$ be the maximal virtually nilpotent subgroup of $\pi$ that 
contains $\pi_e$. 
First, note that any maximal parabolic subgroup $P^\prime$ other than $P$  
fixes a vertex of $T$.
(Indeed, every edge stabilizer of $T$ is a conjugate of $\pi_e$, 
in particular, it lies
in a conjugate of $P$, hence it must have trivial intersection 
with $P^\prime$. So, unless $P^\prime$ fixes a vertex of $T$, 
we get that $P^\prime$ splits over the trivial group
which is impossible because $P^\prime$ is noncyclic virtually nilpotent, and
hence one-ended.)
In particular, if $P$ is not parabolic,
then the $\pi$-action on $T$ satisfies the desired properties.

Now assume that $P$ fixes a vertex $v$ of $T$. Then the fixed-point-set
of $\pi_e$ contains the segment joining $v$ and a vertex of $e$.
In particular, $\pi_e$ lies in $\pi_{\bar e}$, the stabilizer of an edge
$\bar e$ adjacent to $v$. Since $\pi_{\bar e}$ is conjugate to $\pi_e$,
it is virtually nilpotent. Hence $\pi_{\bar e}$ lies in
a maximal parabolic subgroup containing $\pi_e$ which is $P$.
Thus the splitting of $\pi$ over $\pi_{\bar e}$ satisfies the 
desired properties.

It remains to consider the case when 
$P$ is parabolic and $P$ does not fix a vertex of $T$.
Let $\tau$ be the unique $P$-invariant minimal subtree of $T$.
Let $g$ be an element of $\pi$ such that $g^{-1}e$ is an edge of $\tau$
(such a $g$ exists since there is only one orbit of edges).
Then the stabilizer of $g^{-1}e$ in $\pi$ is 
$g^{-1}\pi_e g$. Note that the group $g^{-1}\pi_e g\cap P$ 
(which is the stabilizer of $g^{-1}e$ in $P$) is infinite because
otherwise $P$ splits over a finite group $g^{-1}\pi_e g\cap P$
and $P$ cannot since noncyclic nilpotent groups are one ended.
Since $\pi_e\le P$, we conclude that $g^{-1}Pg\cap P$ is infinite,
therefore $P=g^{-1}Pg$ and, hence, $g\in P$ because
any maximal parabolic subgroup is equal to its normalizer.
The above argument has several implications as follows.

$\bullet$ $e$ must be an edge of $\tau$.
(Indeed, using that there is only one orbit of edges
we find $g\in\pi$ such that $g^{-1}e$ is an edge of $\tau$.
By above $g$ must belong to $P$ and since 
$P$ stabilizes $\tau$,
$e$ is an edge of $\tau$.)

$\bullet$ $P$ is equal to the setwise stabilizer of $\tau$.
(Indeed, if $h\tau =\tau$, than $h$ takes the edge $e$ of $\tau$
to an edge of $\tau$. Hence, again $h\in P$.) 

$\bullet$ Any edge of $T$ lies in a unique $\pi$-image of $\tau$.
(Indeed, any edge of $T$ is $\pi$-equivalent to $e$, so
it is an edge of a tree that is a $\pi$-image of $\tau$. 
Uniqueness of this tree is deduced as follows. 
Suppose, arguing by contradiction,
that two trees that are in $\pi$-image of $\tau$ share an edge. 
Applying an element of $\pi$
we can assume these two trees are $\tau$ and, say, $f\tau$ for some
$f\in \pi$.
This common edge is of the form $g^{-1}e$ for some $g\in\pi$, so
by above $g\in P$.
Hence, $e$ is the common edge of $\tau=g\tau$ and $gf\tau$.
So $e$ is $gf$-image of an edge of $\tau$, or equivalently,
$f^{-1}g^{-1}e$ is an edge of $\tau$. By the same argument, $gf\in P$,
so $f\in P$ which contradicts to $\tau\neq f\tau$.)

Let $Q$ be the set of $\pi$-images of $\tau$. Construct now a 
graph $S(Q)$ with the vertex set $V(T)\cup Q$ 
where we deem two vertices $v,\eta$
adjacent iff $v\in V(T)$, $\eta\in Q$, and $v\in \eta$.

First, prove that $S(Q)$ is a (simplicial) tree. 
Indeed, to show that $S(Q)$ is connected
it suffices to find an arc that joins two arbitrary vertices $x,y\in V(T)$.
The shortest arc that joins $x$ and $y$ in $T$ can be uniquely written  
in the form $x_0x_1\dots x_n$ where $x_0=x$, $x_n=y$ and
the arc $x_{i-1}x_i$ is contained in the tree $\tau_i\in Q$ 
(any edge of $T$ lies in a unique tree from $Q$). 
Then the arc $x_0\tau_1 x_1\dots \tau_n x_n$ joins $x$ and $y$ in $S(Q)$.
Second, show that $S(Q)$ has no circuits.
Indeed, suppose $x_1\tau_1 x_2\dots \tau_n x_n$ is a circuit in $S(Q)$
with $x_1=x_n$. Let $\alpha_i$ be the arc in $\tau_i$ that
connects $x_i$ to $x_{i+1}$. Then $\alpha_1\cup\alpha_2\dots\cup\alpha_n$
is a circuit in $T$, a contradiction.

The group $\pi$ acts on $S(Q)$ without edge inversions. 
Any maximal parabolic
subgroup now fixes a vertex. In particular, $P$ is the stabilizer of $\tau$.
Since $\pi$ is not virtually
nilpotent, there exists a vertex $v$ of $\tau$ whose stabilizer $\pi_v$
is not a subgroup of $P$. (If the stabilizers of two adjacent 
vertex groups of $\tau$ were subgroups of $P$,
then, since $T$ has only one $\pi$-orbit of edges, 
the groups $A$ and $B$ would lie in a conjugate of $P$. 
Hence $\pi=\langle A,B\rangle$ would lie in the conjugate of $P$.)

Now look at the edge $d$ that joins vertices $v$ and $\tau$. 
The stabilizer $\pi_d$ of $d$ is a subgroup of $P$, 
so it is not equal to $\pi_v$.
If $\pi_d=P$, then $P$ has to stabilize
$v$ as well, hence $P$ fixes a point of $T$ which is not the case.
Thus, $\pi_d$ is not equal to the stabilizers of $v$ and
$\tau$. 

Now if any vertex of $T$ is $\pi$-equivalent to $v$, 
then $S(Q)$ has one orbit of edges and
$\pi\cong \pi_v\ast_{\pi_d} P$ is the desired splitting.
Otherwise, $T$ has two orbits of vertices represented by
$v$ and $v^\prime\in\tau$ and hence
$S(Q)$ has two orbits of edges $d$, $d^\prime$.
This defines a splitting 
$\pi\cong \pi_v\ast_{\pi_d} (P\ast_{\pi_{d^\prime}}\pi_{v^\prime})$ 
over $\pi_d$ as needed. 
\end{proof}

\begin{nothing}
\bf{Gluing aspherical cell complexes.}\rm\
In this section we frequently use the following standard construction.
Let $f:X\to Y$ and $g:X\to Z$ be (cellular) maps of cell complexes
that induce $\pi_1$-injections on each connected
component of $X$. Assume that $Y$, $Z$ and each
connected component of $X$ are aspherical.
Form a cell complex by gluing $X\times [0,1]$ to 
$Y$ and $Z$ where $(x,0)$ and $(x,1)$
are identified with $f(x)$ and $g(x)$, respectively. 
The result is an aspherical cell complex 
(use Mayer-Vietoris in the universal covers
and then the Hurewitz theorem). 
When $X$ is connected, its fundamental group 
is isomorphic to $\pi_1(Y)*_{\pi_1(X)} \pi_1(Z)$.
Similarly, if $Y=Z$ we can form a cell complex
by gluing $X\times [0,1]$ to 
$Y$ where $(x,0)$ and $(x,1)$
are identified with $f(x)$ and $g(x)$, respectively.
The result is, again,
an aspherical cell complex and, if $X$ is connected, its
fundamental group is isomorphic to $\pi_1(Y)*_{\pi_1(X)}$.
\end{nothing}

\begin{nothing}
\bf{Topological model for the splitting: amalgamated product case.}\rm\
\label{model amalgam}
Assume that $\pi$ splits over a virtually nilpotent group 
so that~\ref{relative split} gives an isomorphism $A\ast_C B\to\pi$
such that any maximal parabolic subgroup is conjugate
to a subgroup of $A$ or $B$, and the maximal virtually nilpotent 
subgroup $P$ containing $C$ lies in $A$. 
The inclusions of $A$ and $B$ into $\pi$
define coverings $M_A\to M_\pi$, $M_B\to M_\pi$.

Note that for each connected component $L$ of $\partial M_\pi$,
the inclusion $L\hookrightarrow M_\pi$ lifts to $M_A$ or $M_B$.
(Indeed, maps into aspherical manifolds are homotopic iff the induced
$\pi_1$-homomorphisms are conjugate so $L\hookrightarrow M_\pi$
is homotopic to a map that takes
$\pi_1(B)$ to $A$ or $B$. This map lifts to the 
corresponding cover, and hence so does $L\hookrightarrow M_\pi$
by the covering homotopy theorem.)
For every boundary component of $M_\pi$ fix such a lift
thereby defining an injection 
$\pi_0(\partial M_\pi)\to\pi_0(\partial M_A\cup\partial M_B)$.
We refer to these components of $\partial M_A\cup\partial M_B$
as {\it lifted}. 

Consider a closed aspherical manifold $N_P$
with $\pi_1(N_P)\cong P$ defined as follows.
If $P$ is parabolic, $N_P$ is the 
boundary component of $M_\pi$ corresponding 
to the inclusion $P\hookrightarrow\pi$.
If $P$ is loxodromic, $N_P$ is a circle embedded in $M_\pi$
that realizes a generator of $P$. If $P$ is trivial,
$N_P$ is a point of $M_\pi$.

In any case the inclusion $N_P\hookrightarrow M_\pi$ 
lifts to $M_A$, and we fix such a lift (in case $P$ is parabolic,
we use the same lift onto a boundary component of $M_A$
that has been chosen above.)
We use the notation $\tilde N_P$ for 
the lifted manifold.
The covering $N_C\to N_P$ induced by $C\hookrightarrow P$
lifts to a covering onto $\tilde N_P\subset M_A$. 

Also fix a lift of $N_C\to N_P$ to $M_B$. This lift turns out to
be a homeomorphism of $N_C$ onto
a boundary component of $M_B$ because $C=A\cap B=P\cap B$. 
This component is denoted by $\tilde N_C$.

The above maps of $N_C$ into $M_A$ and $M_B$ certainly
become homotopic after projecting to $M_\pi$.
Now we use this homotopy to build an aspherical cell complex $Y$ 
by gluing $M_A$, $M_B$, and
$N_C\times [0,1]$ identifying $N_C\times\{0\}$ with $\tilde N_P$
via the covering $N_C\to\tilde N_P$, 
and $N_C\times\{1\}$ with $\tilde N_C$ via the chosen lift. 
Now we get a homotopy equivalence
$Y\to M_\pi$ which extends the coverings $M_A\to M_\pi$ and  $M_B\to M_\pi$.

Doubling $M_\pi$ along the boundary produces
a closed aspherical manifold $DM_\pi$. Similarly, the
doubles of $M_A$, $M_B$ along the lifted boundary components
are denoted by $DM_A$, $DM_B$.
Finally, let $DY$ be the ``double'' of $Y$ along the 
lifted boundary components, that is,
$DY$ is obtained by gluing $DM_A$, $DM_B$, and two copies of
$N_C\times [0,1]$ as above. The homotopy equivalence
$Y\to M_\pi$ now extends to a homotopy equivalence $DY\to DM_\pi$.

It is useful to give another description of $DY$. 
Namely, let $\bar{D}M_B$ be the manifold obtained from $DM_B$
by identifying two copies of $\tilde N_C\subset\partial DM_B$.
Then $DY$ is obtained by gluing $DM_A$, $\bar{D}M_B$, and 
$N_C\times [0,1]$ where $N_C\times \{1\}$ is identified with
$\tilde N_C\subset\bar{D}M_B$, and $N_C\times \{0\}$ 
is glued to $\tilde N_P\subset M_A$ via the covering 
$N_C\to \tilde N_P$.
\end{nothing}

\begin{nothing}
\bf{Topological model for the splitting: HNN-extension case.}\rm\
As above, fix an isomorphism $A\ast_C\to\pi$
such that any maximal parabolic subgroup is conjugate
to a subgroup of $A$ or $B$, and the maximal virtually nilpotent 
subgroup $P$
containing $C$ lies in $A$. (We think of $A\ast_C$ as
$\langle A,t| \phi(C)=tCt^{-1}\rangle$ where $\phi\co C\to A$
is a monomorphism.)

We keep essentially the same notations as in the amalgamated product case.
Thus the inclusion $A\hookrightarrow\pi$ 
defines a covering $M_A\to M_\pi$ and, 
we fix lifts of all boundary
components of $M_\pi$ to $\partial M_A$.
Again, a lift of $N_P$ to $\partial M_A$ is denoted by $\tilde N_P$.
The covering $N_C\to N_P$ induced by the inclusion 
$C\hookrightarrow P$ of course lifts to 
a covering onto $\tilde N_P$ using the chosen lift of $N_P$.

Also the composition $N_C\to N_P\hookrightarrow M_\pi$ 
can be lifted to a homeomorphism onto a component $\tilde N_P$
of $\partial M_A$ 
so that the lift induces $\phi\co C\to A$.
(Indeed, let $f$ be the lift of $N_C\to N_P\hookrightarrow M_\pi$
to the universal covers such that $f$ is equivariant with respect to
the inclusion $C\hookrightarrow\pi$. Then $g=tf$ is a
$\phi$-equivariant homeomorphism onto a boundary component of the 
universal cover of $M_\pi$. So $g$ descends to a covering
of $N_C$ onto a component of $\partial M_A$ which is, in fact,
a homeomorphism because $\phi(C)$ is maximal parabolic in $A$. 
The last statement is true since
the maximal virtually nilpotent subgroup of $\pi$ containing $\phi(C)$
is $tPt^{-1}$ and $\phi(C)=tAt^{-1}\cap A=tPt^{-1}\cap A$.)

Again,  $\bar{D}M_A$ is a double of $M_A$ along the union of the
lifted boundary components and $\tilde N_C$.
As before we get a homotopy equivalence between $DM_\pi$ and
$DY$ where $DY$ is obtained by gluing $\bar{D}M_A$ and 
$N_C\times [0,1]$ so that $N_C\times \{0\}$ is identified with
a copy of $\tilde N_C$ sitting inside $\bar{D}M_A$, while $N_C\times \{1\}$ 
is glued to $\tilde N_P$ via the covering $N_C\to \tilde N_P$.
\end{nothing}

\begin{thm} \label{nosplit thm}
Let $\pi$ be the fundamental group 
of a finite volume complete Riemannian manifold $M$ of dimension
$n>2$ and with sectional curvatures within $[a,b]$ for $a\le b<0$. 
Then $\pi$ does not split over a virtually nilpotent group.
\end{thm}
\begin{proof} 
Assume first that the splitting of $\pi$ given by~\ref{relative split}
is $A\ast_C B$. Think of $DY$ as glued from $DM_A$, $\bar DM_B$, and
$N_C\times [0,1]$. 
Since the splitting is nontrivial, both $A$ and $B$
have infinite index in $\pi$ so $M_A$, $M_B$ are
noncompact. Hence $DM_A$, $\bar{D}M_B$ are noncompact
since they are glued from noncompact spaces along compact subset.
Now look at the 
Mayer-Vietoris sequence for homology with 
$\mathbb Z/2\mathbb Z$-coefficients.
$$0\to H_n(DY)\to H_{n-1}(N_C)\to  H_{n-1}(DM_A)\oplus
H_{n-1}(\bar{D}M_B)\to\dots$$
The map $H_{n-1}(N_C)\to  H_{n-1}(DM_A)\oplus H_{n-1}(\bar{D}M_B)$
can be written as $i_{a\ast}\oplus -i_{b\ast}$
where $i_a\co N_C\to DM_A$ is the covering 
onto $\tilde N_P$, and
$i_b\co N_C\to \bar{D}M_B$ is the homeomorphism onto $\tilde N_C$.
Show that $i_{b\ast}$ is injective. This is clear if
the cohomological dimension of $C$ is $<n-1$.
Otherwise, by the exact sequence of the pair $(\bar{D}M_B, \tilde N_C)$
it suffices to show that $H_n(\bar{D}M_B, \tilde N_C)=0$ which is 
true by~\cite[VIII.3.4]{Dol}. (As stated,~\cite[VIII.3.4]{Dol}
only applies when the complement of $\tilde N_C$ in $\bar{D}M_B$ is connected.
However, if $\bar{D}M_B\setminus \tilde N_C$ is nonconnected,
it has two noncompact components each equal to $M_B$
and the result again follows
since by the relative Mayer-Vietoris sequence
$H_n(\bar{D}M_B, \tilde N_C)$ is a sum of two copies of 
$H_n(M_B, \tilde N_C)=0$.)

By exactness $H_n(DY)=0$ and we get a contradiction with the fact 
that $DY$ is homotopy equivalent to the closed $n$-manifold $DM_\pi$.

Similarly, if $\pi\cong A\ast_C$, then we get the following 
Mayer-Vietoris sequence with 
$\mathbb Z/2\mathbb Z$-coefficients~\cite{Bro}.
$$0\to H_n(DY)\to H_{n-1}(N_C)\to  H_{n-1}(\bar{D}M_A)\to\dots$$
The map $H_{n-1}(N_C)\to  H_{n-1}(\bar{D}M_A)$ can be written as
$i_*-f_*$
where $i\co N_C\hookrightarrow \bar{D}M_A$ is the homeomorphism 
onto $\tilde N_C$, and
$f$ is the covering of $N_C$ onto $\tilde N_P\hookrightarrow DM_A$. 
It remains to show that $i_*-f_*$
is injective which is clear if the cohomological dimension of $C$ is $<n-1$.
Otherwise, look at 
the exact sequence of the pair $(\bar{D}M_A, \tilde N_C\coprod\tilde N_P)$ 
with $\mathbb Z/2\mathbb Z$-coefficients.
Again, by~\cite[VIII.3.4]{Dol} 
$H_n(\bar{D}M_A, \tilde N_C\coprod \tilde N_P)=0$,
hence the inclusion $\tilde N_C\coprod \tilde N_P\hookrightarrow\bar{D}M_A$ 
induces an injection on $(n-1)$-homology. 
Therefore, $$(i_*-f_*)[N_C]=[N_C]-[N_P]\deg(N_C\to N_P)$$ is nonzero
as wanted an we get a contradiction as before.
\end{proof}

\section{Convergence of finite volume manifolds}
We refer to~\cite{Bel5} or~\cite{Bel7} for background. 
Here we only recall basic definitions and prove several
new lemmas specific to the finite volume case.
By an {\it action} of an abstract group $\pi$ on a space $X$
we mean a group homomorphism $\rho:\pi\to\mathrm{Homeo}(X)$.
An action $\rho$ is called {\it free} if $\rho(\gamma)(x)\neq x$
for all $x\in X$ and all $\gamma\in \pi\setminus\text{id}$.
In particular, if $\rho$ is a free action, then $\rho$ is injective. 

\begin{nothing}
\bf{Equivariant pointed Lipschitz topology.}\rm\
Let $\Gamma_k$ be a discrete
subgroup of the isometry group
of a complete Riemannian manifold $X_k$
and $p_k$ be a point of $X_k$.
The class of all such triples $\{(X_k, p_k, \Gamma_k)\}$ 
can be given the so-called 
equivariant pointed Lipschitz topology~\cite{Fuk1}; 
when $\Gamma_k$ is trivial this reduces to the usual 
pointed Lipschitz topology.

If $(X_k,p_k)$ is a sequence of simply-connected complete
Riemannian $n$-manifolds with $a\le sec(X_k)\le 0$,
then $(X_k,p_k)$ subconverges in the pointed Lipschitz 
topology to $(X,p)$ where $X$ is a $C^\infty$--manifold with 
a complete $C^{1,\alpha}$--Riemannian metric of Alexandrov 
curvature $\ge a$ and $\le 0$.
In fact, in a suitable harmonic atlas on $X$, 
the sequence $(X_k,p_k)$ subconverges to $(X,p)$
in pointed $C^{1,\a}$ topology~\cite{Fuk1, And2}.
\end{nothing}

\begin{nothing}
\bf{Pointwise convergence topology.}\rm\
Suppose that the sequence  $(X_k, p_k)$ converges to $(X,p)$
in the pointed Lipschitz topology.
This allows one to talk about the {\it convergence} of 
a sequence of points $x_k\in X_k$ to $x\in X$
Furthermore, a sequence
of isometries $\gamma_k\in\mathrm{Isom}(X_k)$
we say that $\gamma_k$ {\it converges}, if for any $x\in X$ and any 
sequence $x_k\in X_k$ that converges to $x$,
$\gamma_k(x_k)$ converges. 
The limiting transformation $\gamma$ that takes $x$ to
the limit of $\gamma_k(x_k)$
is necessarily an isometry of $X$.

Let $\rho_k:\pi\to\mathrm{Isom}(X_k)$ 
be a sequence of isometric actions of a group $\pi$ on $X_k$.
We say that a sequence of actions $(X_k, p_k, \rho_k)$ 
{\it converges in the pointwise 
convergence topology} if $\rho_k(\gamma)$ converges 
for every $\gamma\in\pi$.
The map $\rho:\Gamma\to\text{Isom}(X)$ that takes 
$\gamma$ to the limit of $\rho_k(\gamma)$ is necessarily
a homomorphism. 
\end{nothing}

\begin{lem} \label{finite index}
Let $\rho_k:\pi\to\mathrm{Isom}(X_k)$
be a sequence of free isometric actions
of a discrete group $\pi$ on Hadamard 
$n$-manifolds $X_k$ of sectional curvatures within $[a,b]$ for $a\le b<0$.
Assume that the sequence
$(X_k, p_k, \rho_k(\pi))$ converges in the 
equivariant pointed Lipschitz topology to $(X, p, \Gamma)$
and $(X_k, p_k, \rho_k)$ converges to $(X, p, \rho)$ 
in the pointwise convergence topology. 
If $\rho(\pi)$ has finite index in $\Gamma$,
then $\Gamma=\rho(\pi)$.\end{lem}

\begin{proof} 
Being the fundamental group of a pinched negatively curved manifold,
$\pi$ is torsion free and any abelian subgroup of $\pi$
is finitely generated~\cite{BS, Bow2}. 
Therefore, every element in $\pi$ is a power of a primitive element
(an element of a group is called primitive it is not a proper power).
Fix an arbitrary $h\in\Gamma$ and find $g\in\pi$ so that
$h^m=\rho(g)$ for some $m\ge 1$.  
By passing to the appropriate root we can assume $g$ is primitive.
There always exists a sequence $\rho_k(g_k)$ that converges to $h$.
Then both $\rho_k(g_k^m)$ and $\rho_k(g)$ converge to $h^m$.
Hence $g_k^m=g$ for large $k$ by~\cite[lem 2.6]{Bel5}. 
Since $g$ is primitive, $m=1$ so $h\in\rho(\pi)$ as wanted.
\end{proof}

In general, if $|sec|$ is bounded above, the 
injectivity radius is at best upper-semicontinuous.
However, it becomes continuous if the curvature is also nonpositive.

\begin{lem} If $(M_k,p_k)$ be a sequence
of pointed complete Riemannian manifolds of sectional curvatures within
$[a,0]$ for $a\le 0$. Assume that $(M_k,p_k)$ converges to $(M,p)$
in pointed $C^{1,\alpha}$ topology. Then for any
$x_k\in M_k$ that converges to $x\in M$,
the injectivity radii of $M_k$ at $x_k$ converge to
the injectivity radius of $M$ at $x$. 
\end{lem}
\begin{proof} Arguing by contradiction 
find $x_k$ that converges to $x$ while $inj(x_k)$
does not approach $inj(x)$. Pass to a subsequence 
so that $inj(x_k)$ converges to $c\neq inj(x)$. 
Write $M_k=X_k/\Gamma_k$ and $M=X/\Gamma$ and pick
preimages $\tilde x_k$ of $x_k$, and $\tilde x$ of $x$. 
Since $|sec(M_k)|$ is uniformly bounded, a result in~\cite{Fuk1} implies that
$(X_k, \tilde x_k, \Gamma_k)$ subconverges to
$(X,\tilde x,\Gamma)$ in the equivariant pointed
Lipschitz topology.
Now since the sectional curvature is nonpositive,
$inj(x)=d_\Gamma(\tilde x)/2$ and $inj(x_k)=d_{\Gamma_k}(\tilde x_k)/2$ 
where $d_\Gamma(x)$ is the infimal displacement of the point
$x$ by the isometries of $\Gamma$. 
It follows easily from the definition of equivariant
Lipschitz convergence that $d_{\Gamma_k}(\tilde x_k)$ converges to 
$d_\Gamma(x)$ so we get a contradiction.
\end{proof}

\begin{nothing}
\bf{Nikolaev's smoothing theorem.}\rm\
A useful technical tool is the following smoothing
theorem of I.~Nikolaev~\cite{Nik}. We only state an easy 
case of $C^{1,\a}$ metrics even though~\cite{Nik}  
actually applies to any path metric.
Let $(M, g)$ be a complete $C^{1,\a}$-Riemannian
manifold of Alexandrov curvature $\ge a$ and $\le b$.
Then there exist a sequence of complete
Riemannian metrics $g_m$ on $M$ with 
sectional curvatures within $[a-1/m,b+1/m]$ such that
$(M,g_m)$ converges to $(M,g)$ in (unpointed!)
$C^{1,\a}$-topology. Also
$id:(M,g_m)\to (M,g)$  and $id:(M,g)\to (M,g_m)$ 
are $2^{1/m}$-Lipschitz. 
\end{nothing}

\begin{prop}\label{convergence lemma}
Let $X_k$ be a sequence of Hadamard $n$-manifolds
with sectional curvatures in $[a, b]$ for $a\le b<0$ and $n\ge 3$.
Let $\rho_k: \pi\to\mathrm{Isom}(X_k)$ be  
an arbitrary sequence of free and isometric actions
such that each $X_k/\rho_k(\pi)$ has finite volume.
Then after passing to a subsequence, there are points $p_k\in X_k$ such that 
\begin{description}
\item{(i)} $(X_k, p_k, \rho_k)$ converges to 
in the pointwise convergence topology to a free action $(X, p, \rho )$, and
$(X_k, p_k, \rho_k(\pi))$ converges
in the equivariant pointed Lipschitz topology to  $(X, p, \rho (\pi))$;
\item{(ii)} If $q_k$, $q$ are projections of $p_k$, $p$, respectively,
then $(X_k/\rho_k(\pi), q_k)$ converges to $(X/\rho(\pi),q)$
in the pointed $C^{1,\a}$ topology;
\item{(iii)} 
the $C^{1,\a}$-manifold $X/\rho(\pi)$ has finite volume, and 
is diffeomorphic to the
interior of a compact manifold with boundary; the number
of connected components of the boundary is the number
of maximal parabolic subgroups of $\pi$.
\end{description}
\end{prop}
\begin{proof}
Let $M_k=X_k/\rho_k(\pi)$.
By~\ref{nosplit thm} and~\cite[2.7, 2.10]{Bel5} 
we can assume by passing to a subsequence that
$(X_k, p_k, \rho_k)$ converges to $(X, p, \rho)$ in pointwise
convergence topology and $(X_k, p_k, \rho_k(\pi))$ converges to 
$(X, p, \Gamma )$ in equivariant pointed Lipschitz
topology. By~\cite[2.5]{Bel5} we know that the $\Gamma$-action on $X$
is free.
By~\cite{Fuk1, And2} $(M_k,q_k)$ converges to $(X/\Gamma,q)$
in the pointed $C^{1,\a}$ topology.

Now $M=X/\rho(\pi)$ is a smooth manifold with 
a $C^{1,\alpha}$ Riemannian metric $g$
of Alexandrov curvature of $g$ is within $[a,b]$.
Let $(M,g_m)$ be the Nikolaev's smoothing of $(M,g)$~\cite{Nik}.
Since $(M,g_m)$ is a complete Riemannian manifold of
pinched negative curvature which is homotopy equivalent
to $M_k$,~\cite{Sch} implies that $(M,g_m)$ has finite volume
so~\cite{BGS} implies that $M$ is the interior of a compact manifold
with boundary and the number of boundary components
is the same as the number of ends of $M_k$ (or, alternatively,
the number of maximal parabolic subgroups of $\pi$).

Since  $(M,g_m)$ converge to $(M,g)$ in unpointed
$C^{1,\a}$-topology, $vol(M,g_m)$ tends to $vol(M,g)$.
Also $id:(M,g_m)\to (M,g_{m+1})$ is Lipschitz with Lipschitz
constant approaching $1$, hence
the volumes of $(M,g_m)$ are uniformly bounded~\cite[p11]{BGS}.
Thus the volume of $(M,g)$ is finite.
Being a quotient of $M$, the manifold $X/\Gamma$
also has finite volume so $\rho(\pi)$ has 
finite index in $\Gamma$.
Hence $\rho(\pi)=\Gamma$ by~\ref{finite index}.
\end{proof}

\begin{rmk} \label{nikolaev and thin part}
Another obvious application of Nikolaev's smoothing and the continuity of
injectivity radius is that
$$\{inj_{g_k}\le\epsilon 2^{-1/k}\}\subset
\{inj_g\le\epsilon\}\subset
\{inj_{g_k}\le\epsilon 2^{1/k}\}.$$
In particular, $\epsilon$-thick part 
$\{inj_{g}\ge\epsilon\}$ is compact
since it is so for $C^\infty$ metrics.
\end{rmk}

\begin{lem} \label{diam of thick part}
The injectivity radius and the
diameter of $\epsilon$-thick part are
bounded above and below on $\mathcal M_{a, b,\pi, n}$.
\end{lem}
\begin{proof} 
Arguing by contradiction, let $M_k\in \mathcal M_{a, b,\pi, n}$ 
be a sequence of manifolds with 
the diameter of the $\epsilon$-thick part going to 
infinity for some $\epsilon$.
Find points $p_k\in M_k$ such that $(M_k, p_k)$ subconverges
to $(M,p)$. 
By~\ref{nikolaev and thin part}
the $\epsilon$-thick part of $M$ is compact hence it lies 
in the open ball $B(p,R)$ for some $R$. 
Pass to subsequence so that diameters of 
the $\epsilon$-thick parts of $M_k$ are $>2R+2$; thus, for all large
$k$, the $\epsilon$-thick part
of $M_k$ contains a point that lies in $M_k\setminus B(p_k,R+1)$.
Furthermore, by continuity of injectivity radius
$B(p_k,R+1)$ contains a point of injectivity radius 
$\ge\epsilon$ for all large $k$.
Since $n>2$, the $\epsilon$-thick part is connected,
so  for all large $k$, 
$B(p_k,2R+1)\setminus B(p_k,R+1)$ contains a point of injectivity radius
$\ge\epsilon$. 
This point subconverges to a point in $B(p,2R+2)\setminus B(p,R+1/2)$
of injectivity radius $\ge\epsilon$ which is a contradiction.

Now a lower bound on the injectivity radius
(and hence on the diameter of the $\e$-thick part) 
is provided by the Margulis lemma~\cite{BGS}.
As for an upper bound assume there is a sequence of manifolds $(M_k,p_k)$ in 
with $M_k\in \mathcal M_{a, b,\pi, n}$ and $inj_{p_k}>k$. 
Pass to a subsequence so that $(M_k,p_k)$ converges to $(M,p)$.
Since the $\epsilon$-thick part of $M$ is compact,
the injectivity radius of $M$ is bounded above
which contradicts the continuity of the injectivity radius. 
\end{proof}

\begin{lem}\label{consists of cusps}
There exist $\e>0$
such that for any  $M\in \mathcal M_{a, b,\pi, n}$
the $\e$-thin part of $M$ consists only of cusps.
\end{lem}
\begin{proof}
Arguing by contradiction find a sequence of manifolds
$M_k=X_k/\rho_k(\pi)\in \mathcal M_{a, b,\pi, n}$
each containing a closed geodesic of length $<1/k$. 
By~\ref{convergence lemma} and~\cite{Fuk1}
$M_k$ subconverges in the pointed
Lipschitz topology. Denote the limiting manifold by $(M,g)=X/\rho(\pi)$.

By~\ref{nikolaev and thin part} the $\e$-thick part of
$(M,g)$ contains the $2\e$-thick part of $(M,g_1)$.
For large $k$ Lipschitz approximations $\phi_k$ between $M_k$ and $M$ take
the $\e$-thick part of $M_k$ Hausdorff close to the $\e$-thick part of $M$.
Hence for all large $k$ the $2\e$-thick part of $(M,g_1)$ 
is contained in the $\phi_k$-image of the $\e$-thick part of $M_k$.

Now let $\e$ be so small that $2\e $-thin part of $(M,g_1)$ 
consists only of cusps. 
Thus, the $\phi_k$-image of any $\e$-tube of $M_k$ lies in a 
$2\e$-cusp of $(M,g_1)$.
Hence if $\gamma_k$ be a closed geodesic in a tube of $M_k$, 
then $\phi_k(\gamma_k)$ represents a parabolic element.
By algebraic reasons, $n>2$ implies that $\rho_k\circ\rho^{-1}$
takes parabolics to parabolics. So $\gamma_k$ represents
a parabolic element which is a contradiction. 
\end{proof}

\section{Group theoretic applications}

In this section $\pi$ is a group such that 
$\mathcal M_{a, b, \pi, n}$ is nonempty.

\begin{cor}
$Out(\pi )$ is finite.
\end{cor}
\begin{proof}
Let $\rho_k\in\mathrm{Aut}(\pi)$
lie in different conjugacy classes. This defines a sequence $\rho_k$
of free isometric actions of $\pi$ on the universal cover $X$ of $M$.
For a finite generating set $S$ of $\pi$, let 
$D_k(x)=\max_{\gamma\in S}\{d(x,\rho_k(\gamma)(x))\}$ and
$D_k=\inf\{D_k(x): x\in X\}$.
As in~\cite[2.10]{Bel5}, we can choose a sequence of points
$x_k\in X$ such that $D_k(x_k)\le D_k+1/k$. 
If $D_k\to\infty$, we get an action on an $\mathbf R$-tree hence a 
splitting which is impossible. So assume $D_k(x_k)$ is 
uniformly bounded.

Let $F$ the the Dirichlet fundamental domain 
for the action of $\pi$ on $X$. 
There exists a sequence $\phi_k\in\pi$ such that
$\phi_k(x_k)\in F$. 
Then $(X, \phi_k(x_k), \phi_k\rho_k\phi_k^{-1})$ subconverges
in the pointwise convergence topology.

Suppose first that
$\phi_k(x_k)$ is precompact, so that passing to subsequence
we can assume that $\phi_k(x_k)$ converges to $x\in X$.
Then $(X, x, \phi_k\rho_k\phi_k^{-1})$ also subconverges
in the pointwise convergence topology, in other words,
passing to subsequence, we deduce that 
$\phi_k\rho_k(\gamma)\phi_k^{-1}$ converges in $\mathrm{Isom}(X)$
for any $\gamma\in\pi$. 
Since $\pi$ is a closed subgroup the limit lies in $\pi$.
So $\phi_k\rho_k\phi_k^{-1}$ converges in 
$\mathrm{Hom}(\pi ,\pi )$. 
Since the space $\mathrm{Hom}(\pi ,\pi )$ is discrete,
$\phi_k\rho_k\phi_k^{-1}$ are all equal for large $k$.
So $\rho_k$ lie in the same conjugacy class for large $k$, a contradiction.

If $\phi_k(x_k)$ is not precompact, then passing to
subsequence we can assume $\phi_k(x_k)$ converges a parabolic
fixed point. (The closure of $F$ at infinity is just finitely
many parabolic fixed points.) Choose a $\pi$-invariant set of
mutually disjoint horoballs. Passing to subsequence,
we assume that $\phi_k(x_k)$ lies in one horoball $H$
for all $k$. Since $\pi$ is not virtually nilpotent,
for each $k$ there is a generator $\gamma_k\in S$ such that
the horoballs $H$ and $\phi_k\rho_k\phi_k^{-1}(\gamma_k )(H)$
are disjoint. 
As $S$ is finite we can pass to to subsequence
so that $H$ and $\phi_k\rho_k\phi_k^{-1}(\gamma_k )(H)$ are disjoint
for all $k$ and some $\gamma\in S$.
Since $D_k(x_k)$ is uniformly bounded, 
the distance between $\phi_k(x_k)$ and 
$\phi_k\rho_k(\gamma )\phi_k^{-1}(\phi_k(x_k))$ 
is uniformly bounded. 
On the other hand, this distance 
has to converge to infinity because it is bounded below by
the distance from to $\phi_k(x_k)$ to the horosphere 
$\partial H$. This contradiction completes the proof.
\end{proof}

\begin{cor}
$\pi $ is cohopfian.
\end{cor}

\begin{proof} Arguing by contradiction assume there exists
an injective homomorphism $\phi:\pi\to\pi$ which is not onto.
By~\cite{Sch} the manifold $X/\phi(\pi)$ has finite volume
hence $\phi(\pi)$ must be of finite index, say $m>1$, in $\pi$.   
Iterating $\phi$, we get a sequence
of free isometric actions $\rho_k$ of $\pi$ on the universal cover $X$ of $M$
such that $\rho_{k+1}(\pi )$ is an index $m$ subgroup of 
$\rho_k(\pi )$ for each $k$.

Same proof as before gives that  
$\phi_k\rho_k\phi_k^{-1}$ are all equal for large $k$. 
In particular, $X/\rho_k(\pi)$ and $X/\rho_{k+1}(\pi)$ 
are isometric for all large $k$.
But there exists an $m$-sheeted cover $X/\rho_{k+1}(\pi)\to X/\rho_k(\pi)$,
hence $vol(X/\rho_{k+1}(\pi))=m\cdot vol(X/\rho_k(\pi))$.
Thus, $m=1$, a contradiction. 
\end{proof}

\begin{rmk} There is of course another proof that $\pi$ is cohopfian. 
Namely, by~\cite{Sch} $\phi(\pi)$ has finite index $m$ in $\pi$.
Since $X/\pi$ and $X/\phi(\pi)$ are properly homotopy equivalent,
they have equal simplicial volumes (which are also nonzero~\cite{Gro5}).
On the other hand, simplicial volume is multiplicative under
finite covers so $m=1$.
\end{rmk}

\section{Diffeomorphism finiteness}

In this section we prove a theorem that implies~\ref{intro: main thm}
when combined with~\ref{convergence lemma}. 
Main ingredients of the proof are exponential convergence of geodesics 
and continuity of injectivity radius. We also use in a crucial
way that the metrics converge in at least $C^1$ topology. 

First,
we need a better understanding of cusps 
for manifolds in $\mathcal M_{a, b, \pi, n}$. 
Fix an $\e\in (0,\mu_{n,a})$ where $\mu_{n,a}$ is the Margulis constant
and fix an $\e$-cusp of a manifold $M=X/\pi\in\mathcal M_{a, b,\pi, n}$.
Denote by $inj_\e$ the boundary of the cusp; 
thus $inj_\e$ is a compact topological submanifold of $M$
of codimension one.
Set $d_\e=\mathrm{diam}(inj_\e)$. Also let $H_\e^+$, $H_\e^-$
be the quotients of horospheres such that 
$inj_\e$ lies in the region bounded by $H_\e^+$ and $H_\e^-$, 
and $H_\e^+\cap inj_\e$, $H_\e^-\cap inj_\e$ are nonempty.
We also assume that $H_\e^-$ is ``closer'' to infinity than $H_\e^+$
(i.e. $\mathrm{diam}(H_\e^-)\le\mathrm{diam}(H_\e^+)$).
Following~\ref{prelim: exp conv} we let $a=-K^2$, $b=-k^2$.

\begin{lem} \label{cusp lemma}
Given $\e,\s$ satisfying
$0<\s<\e<\min\{\mu_{n,a},1\}$, fix an $\e$-cusp of a 
manifold $M=X/\pi\in\mathcal M_{a, b,\pi, n}$. 
Then the following holds:

(1) $\mathrm{diam}(H_\e^+)\le 3d_\e$;

(2) $\mathrm{diam}(inj_\s)\le 2d_\e+k^{-1}\ln\left(\frac{2\e}{\s}\right)-
K^{-1}\ln\left(\frac{\e}{2\s}\right)$ 

(3) $K^{-1}\ln\left(\frac{\e}{2\s}\right)-d_\e\le\mathrm{dist}(inj_\e,inj_\s)\le
k^{-1}\ln\left(\frac{2\e}{\s}\right)+2d_\e$;

(4) if $\s<\frac{\e}{2} e^{-Kd_\e}$, then
$\mathrm{diam}(H_\s^+)\le\mathrm{diam}(H_\e^-)$.
\end{lem}
\begin{proof}
To prove (1) find $x,y\in H_\e^+$ with $dist(x,y)=\mathrm{diam}(H_\e^+)$.
Let $\tilde x, \tilde y\in inj_\e$ be the points obtained by pushing
$x,y$ along radial geodesics. Then by the triangle inequality
$dist(x,y)\le dist(x,\tilde x)+dist(\tilde x,\tilde y)+dist(\tilde y,y)\le 3d_\e$
where the latter inequality holds because Busemann functions are $1$-Lipschitz. 

It remains to prove $(2)-(4)$. Let $z$ be a point of the ideal boundary of $X$
corresponding to the cusp under consideration.
We denote by $\widetilde{inj}_\d$, $\widetilde{H}_\d^+$, $\widetilde{H}_\d^-$
the lifts of $inj_\d$, $H_\d^+$, $H_\d^-$
to the universal cover which are in bounded distance from
a horosphere about $z$. 
Let $\g (t)$ be a geodesic asymptotic to $z$ with $\g(0)\in\widetilde{H}_\e^+$
and assume that $\g(t)$ intersects
$\widetilde{inj}_\s$ and $\widetilde{inj}_\e$ in the points $\g(t_\s)$ and $\g(t_\e)$, 
respectively. Note that $t_\s>t_\e$ and $t_e\in [0,d_\e]$.

Since $inj(\g(t_\s))=\s$, one can find $g\in\pi$ such that 
$d(g(\g(t_\s)), \g(t_\s))=2\s$.
Now $d(g(\g(t_\e)),\g(t_\e))\ge 2\e$. Thus, by~\ref{prelim: exp conv}
\[\frac{\e}{\s}=\frac{2\e}{2\s}\le 
\frac{d(g(\g(t_\e)),\g(t_\e))}{d(g(\g(t_\s)), \g(t_\s))}\le 2e^{K(t_\s-t_\e)}
\ \ \mathrm{or}\ \ K^{-1}\ln\left(\frac{\e}{2\s}\right)\le t_\s-t_\e.\]

Similarly, since $inj(x_\e)=\e$,
one can find $h\in\pi$ with $d(h(\g(t_\e)),\g(t_\e))=2\e$.
Again, $2\s\le d(h(\g(t_\s)),\g(t_\s))$ so
\[\frac{\e}{\s}\ge \frac{d(h(\g(t_\e)),\g(t_\e))}{d(h(\g(t_\s)),\g(t_\s))}\ge 
e^{k(t_\s-t_\e)}/2
\ \ \mathrm{or}\ \ k^{-1}\ln\left(\frac{2\e}{\s}\right)\ge t_\s-t_\e.\]
So $t_\s\in [K^{-1}\ln\left(\frac{\e}{2\s}\right),
d_\e+k^{-1}\ln\left(\frac{2\e}{\s}\right)]$. 
Hence,
by the triangle inequality
\[\mathrm{diam}(inj_\s)\le 2d_\e+k^{-1}\ln\left(\frac{2\e}{\s}\right)-
K^{-1}\ln\left(\frac{\e}{2\s}\right)
\ \ \mathrm{and}\]
\[K^{-1}\ln\left(\frac{\e}{2\s}\right)-d_\e\le\mathrm{dist}(inj_\e,inj_\s)\le
2d_\e+k^{-1}\ln\left(\frac{2\e}{\s}\right).\]
Furthermore, if $K^{-1}\ln\left(\frac{\e}{2\s}\right)>d_\e$, then
$H_\s^+$ is ``closer'' to infinity than $H_\e^-$ as desired.
\end{proof}

\begin{cor}
The volume function is uniformly bounded 
above on $\mathcal M_{a, b,\pi, n}$.
\end{cor}
\begin{proof}
First, show that the volume of the $\s$-thick part is uniformly bounded above
on $\mathcal M_{a, b,\pi, n}$ for any $\s\in (0,\mu_{a,n})$.
Indeed, observe that the diameter of $M_{[\s,\infty)}$ is bounded 
above by~\ref{diam of thick part}.
Hence, $M_{[\s,\infty)}$ is in the image of a ball in $X$ of some
uniformly bounded above radius. By Bishop-Gromov volume comparison
the volume of the ball is uniformly bounded above and the result follows because
the projection $X\to M$ is volume non-increasing.

We fix an $\e\in (0,\mu_{a,n})$, and an 
arbitrary $\e$-cusp of $M\in \mathcal M_{a, b,\pi, n}$. 
Now we seek to obtain a uniform upper bound on $vol(H_\e^-)$.
Let $\s=\frac{\e}{2} e^{-K(d_\e+1)}$ so that 
$H_\s^+$ is ``closer'' to infinity than $H_\e^-$.
Let $T$ be the distance between between $H_\s^+$ and $H_\e^-$;
note that $T\in [d_\e+1, (d_\e+1)K/k+(\ln{4})/k]$.
By above, the volume
enclosed between $inj_\e$ and $inj_\s$ is bounded above on
$\mathcal M_{a, b,\pi, n}$ by a constant $V$ depending only on
$\s,a,b,\pi, n$.
The same is then true for the volume enclosed between
$H_\s^+$ and $H_\e^-$ which is equal to
$\int_{0}^{T} vol(H_t)dt$ where $H_t$ is the quotient of a horosphere 
at $t$-level, and $H_0=H_\e^-$, $H_T=H_\s^+$.

By~\ref{prelim: exp conv}, pushing along radial vector field 
gives an $e^{Kt}$-Lipschitz map $H_t\to H_0$.
Thus, $vol(H_0)\le vol(H_t)e^{Knt}$ 
(in this proof we always equip $H_t$ with the Riemannian metric
induced by the inclusion into $M$). 
Hence by the Fubini's theorem
(which applies since the Busemann function is a $C^2$-Riemannian submersion)
we have
\[vol(H_0)T=\int_{0}^{T}vol(H_0)dt\le e^{KnT}\int_{0}^{T}vol(H_t)dt\le
e^{KnT}V.\]
Thus $vol(H_0)=vol(H_\e^-)$ is uniformly bounded 
above over $\mathcal M_{a, b,\pi, n}$.

Now pushing along radial vector field 
gives an $e^{-kt}$-Lipschitz diffeomorphism $H_0\to H_t$ so
$vol(H_t)\le vol(H_0)e^{-ktn}$. Hence
\[\int_{0}^{\infty}vol(H_t)dt\le vol(H_0)\int_{0}^{\infty}e^{-ktn}dt\le
vol(H_0)/kn\]
and we get a uniform upper bound for the volume
of the $\s$-cusp under consideration. Thus, we get a uniform upper bound
on the volume of $M$.
\end{proof}

\begin{thm} 
Let $(M_k,g_k)$ be a sequence of manifolds in $\mathcal M_{a, b,\pi, n}$
such that for some $p_k\in M_k$,
$(M_k,p_k)$ converges in the pointed $C^{1,\a}$-topology
to $(M, p)$ equipped with a $C^{1,\alpha}$ Riemannian metric $g$.
Then for all large $k$ there exist diffeomorphisms $f_k\co M\to M_k$ 
such that the pullback metrics $f_k^\ast g_k$ 
converge to $g$ in $C^{1,\alpha}$ topology
uniformly on compact subsets. 
\end{thm}
\begin{proof}
Renumerate the sequence so that $M_k$ now come with even indices
while odd indices correspond to Nikolaev's smoothings of $(M,g)$.
We still denote the obtained sequence by $M_k$.
Since $(M_k,p_k)$ converges to $(M,p)$ in the pointed $C^{1,\a}$-topology
there are $1/k$-Lipschitz
smooth embeddings $\phi_k\co B_{1/k}(p_k)\to M$
with $d(p,\phi_k(p_k))\le 1/k$ such that
$\phi_k$-pushforward of $g_k$ converges to $g$
in $C^{1,\a}$-topology uniformly on compact subsets. 
By~\cite{Nik} we can take 
$\phi_{2k+1}=\mathrm{id}_M$.

As we shall prove below, for each small enough $\e$ and large enough $k$, 
there exists a diffeomorphism
$h_{k,\e}\co M_k\to M$ that is equal to 
$\phi_k$ when restricted to the $\e$-thick part of $M_k$.
Note that the continuity of injectivity radius implies that,
given $\e$,  the maps $\phi_k$ are defined on the
$1$-neighborhood of the $\e$-thick part of $M_k$ for all $k\ge C(\e)$
and some positive integer valued function $C$.
Now the diffeomorphism $f_k=h_{k+C(1/k),1/k}$ enjoys the desired properties.
Thus, it remains to construct $h_{k,\e}$.

Use~\ref{consists of cusps} to find an $\e^\prime$ such that 
$\e^\prime$-thin part of each $M_k$ consists of cusps; 
we assume $\e\in (0,\max\{1,\e^\prime/10\})$. 
Fix a cusp of $M$ and the corresponding cusps of $M_k$.
Using~\ref{cusp lemma}, we can make so $\e$ is so small that 

$\bullet$ $H_{k,\e^2}^+$ is closer to infinity than $H_{k,\e}^-$,
and $dist(H_{k,\e^2}^+,H_{k,\e}^-)>10$;

$\bullet$ $H_{k,\e^3}^+$ is closer to infinity than $H_{k,\e^2}^-$, and
and $dist(H_{k,\e^3}^+,H_{k,\e^2}^-)>10$.

(Here subindex $k$ indicates that the quotients of the horospheres 
lie in a cusp of $M_k$.)
Let $R_k$ be the radial
vector field defined on the $\e$-thin part of $M_k$, and 
let $D_k$ be the $1$-neighborhood of the region between 
$inj_{k,\e^2}$ and $inj_{k,\e^3}$. 

Using exponential convergence of geodesics, one can easily see that for any
$\a\in (0,\pi/2)$ there exists $r$, depending only on $\a, a, b, n$ and
independent of $k$, such that for any $x\in D_k$ and any
$y$ lying in the same $\e^2$-cusp of $M_k$ with $d_k(x,y)\ge r$,
the angle at $x$ formed by $R_k$ and the tangent vector to
the geodesic segment $[y,x]$ is $\le\a$.
Fix $\a=\pi/3$ and fix the corresponding $r$.

Assume $k$ is large enough so that the embeddings $\phi_k$
are defined on the $2r$-neighborhood of $D_k$.
By continuity of injectivity radius the domains 
$\phi_k(D_k)\subset M$ converge to some compact set
in Hausdorff topology and we can find a smooth domain $D$
which is Hausdorff close to the set and satisfies 
$D\subset \phi_k(D_k)$ for all large $k$. 
Now use $\psi_k=\phi_k^{-1}$ to pullback all metrics to $D$.

We want to show that the region bounded by $\phi_k(H_{k,\e^2}^-)$ and 
$\phi_{k+1}(H_{k+1,\e^3}^+)$ is diffeomorphic to
$H_{k,\e^2}^-\times [0,1]$. 
It suffices to produce a $C^\infty$
nowhere vanishing vector field on $D$ that is transverse
to both  $\phi_k(H_{k,\e^2}^-)$ and 
$\phi_{k+1}(H_{k+1,\e^3}^+)$. 
We shall construct such a vector field as a controlled approximation 
of $\psi_k^\# R_k$, the $\psi_k$-pullback of the radial vector field $R_k$.
Choose a harmonic atlas $\{B_j\}$ on $D$ 
as in~\cite{And2} in which the $\psi_k$-pullbacks 
of $g_k$ converge to $g$ in $C^{1,\a}$-topology. 
Fix a partition of unity associated with the atlas.
Let $y\in M$ be a point in the same cusp
and $d(y,D)\ge r+1$. 

Now we construct a $C^\infty$ vector field $X_k$ on $D$ 
by defining it on each chart neighborhood $B_j$
and then gluing via the partition of unity.
Look at the $2r$-neighborhood of $D$ equipped with
the metric $\psi_k^\# g_k$.
The preimage of $B_j$ under 
$exp_{\psi_k^\# g_k}\co T_y M\to M$ is the disjoint union 
of copies of $B_j$. Pick a copy closest to the
origin, join each of its points to the origin by rays, and then
project the rays to $M$ via $exp_{\psi_k^\# g_k}$. Now the tangent vectors
at the endpoints of the obtained geodesic segments 
joining $y$ with points of $B_j$ form a vector field
on $B_j$. Gluing these local data via the partition of unity gives $X_k$.
Note that by construction $X_k$ is a nowhere vanishing vector field
such that the angle formed by $X_k$ and the exterior normal $\psi_k^\# R_k$
to $\phi_k(H_{k,\e^2}^-)$ or $\phi_k(H_{k,\e^3}^+)$ is within $(0, \pi/3]$.

Now on each chart $X_k$ is a solution of the geodesic equation.
Since metrics converge in $C^{1,\a}$ topology,  Christoffel symbols
converge in at least $C^0$ topology, so by standard ODE results~\cite[I.5.8]{Reid}
$X_k$ converges in $C^0$-topology to some $C^0$ vector field $X$.  
So the angle, measured in the metric $\psi_k^\# g_k$, 
formed by $X_m$ and $\psi_k^\# R_k$ is within $(0, c]$ for all $m, k$ large enough,
and some $c\in [\pi/3,\pi/2)$.
Now a standard differential topology arguments implies that 
the region between $\phi_k(H_{k,\e^2}^-)$ and 
$\phi_{k+1}(H_{k+1,\e^3}^+)$ is diffeomorphic to
$H_{k,\e^2}^-\times [0,1]$ as needed. 

Finally, we are ready to define $h_{k,\e}$.
Let $M_{k,\d}$ be the compact manifold obtained from $M_k$
by chopping off cusps along all surfaces $H_{k,\d}^-$;
we think of $M_{k,\d}$ as a bounded domain in $M_k$. 
Define $h_{2k+1,\e}=\mathrm{id}_M$, and
define $h_{2k,\e}$ as the following composition.
First, map $M_{2k}$ diffeomorphically to the interior of
$M_{2k,\e^2}$ by a map which is the identity on the $1$-neighborhood
of $M_{2k,\e}$. Then map $M_{2k,\e^2}$ to $M$ by $\phi_{2k}$.
Next use the above argument to map $\phi_{2k}(M_{2k,\e^2})$
diffeomorphically onto $\phi_{2k+1}(M_{2k+1,\e^3})=M_{2k+1,\e^3}$
by a map which is the identity on the $1$-neighborhood
of $M_{2k+1,\e}$. 
Last, map $M_{2k+1,\e^3}$ diffeomorphically to $M_{2k+1}=M$,
again keeping $M_{2k+1,\e}$ fixed.
\end{proof}

\section{Pinching}

In this section we prove several pinching results that follow 
from~\ref{intro: main thm}. 

\begin{cor}\label{pinch: hyperbolic pinching}
Given a group $\pi$, there exists 
$\epsilon=\epsilon(\pi) >0$
such that any finite volume manifold from $\mathcal M_{n,-1-\epsilon,-1,\pi}$
is diffeomorphic to a real hyperbolic manifold.
\end{cor}
\begin{proof}
Arguing by contradiction find a sequence $M_k$ of
manifolds with fundamental group isomorphic to
$\pi$ and sectional curvatures within $[-1-1/k,-1]$.
By the main theorem, we can assume that $(M_k, g_k)$ converges in 
$C^{1,\alpha}$ topology to a $C^{1,\alpha}$-Riemannian manifold $(M,g)$. 
Now the universal cover $X$ of $M$ is a complete $C^{1,\alpha}$-Riemannian
manifold of Alexandrov curvature $-1$. By~\cite{Ale}
$X$ is isometric to the real hyperbolic space and we are done.
(Alternatively, one can repeat the argument below appealing to~\cite{Gao}
to deduce that $g$ is a $C^\infty$ metric.)
\end{proof}

\begin{rmk}
It follows from~\cite{Gro5} that $\e$ in~\ref{pinch: hyperbolic pinching} 
only depends on the simplicial
volume of $M$ when $n>3$. This is formally stronger that dependence on $\pi_1(M)$,
yet it actually amounts to the same thing because the bounds 
$||M||<V$, $a\le sec(M)\le b<0$ imply finitely 
many possibilities for $\pi_1(M)$~\cite{Fuk2}. 
Furthermore,
if $n$ is even and $>3$, then $\e$ in~\ref{pinch: hyperbolic pinching} only depends
on the Euler characteristic $\chi(M)$. Indeed, it is well known that
there is a positive constant $C_n$ such that if $-1-C_n\le sec(M)\le -1$, then
the Gauss-Bonnet integrand $\chi$ satisfies $C_{n,1}\omega\le\chi\le C_{n,2}\omega$
for some constants $C_{n,1}$ ,$C_{n,2}$ and the Riemannian volume form $\omega$. 
Now the Gauss-Bonnet formula 
(generalized to finite volume case in~\cite{CG})
shows that bounds on $vol(M)$ and $\chi(M)$ are equivalent.
\end{rmk}

\begin{rmk}
Note that if $n=4$, then the Gauss-Bonnet formula 
(again see~\cite{CG} in noncompact case) gives a stronger 
version of~\ref{intro: main thm}. Namely, the class
of finite volume Riemannian $4$-manifolds with sectional curvatures
within $[a,b]$ where $a\le b<0$ and with uniformly bounded
Euler characteristics is precompact in the $C^{1,\a}$-topology. 
Indeed, it was shown by Milnor~\cite{Che} that in dimension four 
the Gauss-Bonnet integrand $\chi$
satisfies $\chi\ge 3b^2\omega$ where $\omega$ is the volume form
and $b<0$ is the upper curvature bound.
Hence the Gauss-Bonnet formula implies that, 
for any closed Riemannian $4$-manifold $M$ with sectional curvatures
within $[a,b]$ where $b<0$, volume is bounded by
the Euler characteristic: $\chi(M)\ge 3b^2vol(M)$.
\end{rmk}

Berger's classification of holonomy groups implies the following
possibilities for the restricted holonomy group of a complete negatively
curved $n$-manifold: $\mathbf{SO}(n)$ (generic case), 
$\mathbf{U}(n/2)$ (K\"ahler), $\mathbf{Sp}(1)\mathbf{Sp}(n/4)$
(quaternionic-K\"ahler), and $\mathbf{Spin}(9)$ (Cayley).
Any complete negatively
curved manifold with restricted holonomy group $\mathbf{Spin}(9)$
is a quotient of the Cayley hyperbolic plane~\cite[10.96.VI]{Bes}.
The following uniformization theorem is due to
Yeung~\cite{Yeu}. Let $M$ be a complete finite volume Riemannian manifold
that is either K\"ahler or quaternionic-K\"ahler.
In case $M$ is noncompact assume also that sectional
curvatures of $M$ are within $[a,b]$ for some $a\le b<0$
with $a/b\le (n-1)^4$. If $M$ is homotopy equivalent to
a complete pointwise quarter-pinched negatively curved Riemannian 
manifold, then $M$ is locally symmetric.

Combining~\ref{intro: main thm}, the Yeung's theorem, 
some results of Gao~\cite{Gao}, 
and some linear algebra of K\"ahler
curvature tensor, we get the following 

\begin{cor} \label{kahler pinching}
Given a group $\pi$, there exists 
$\epsilon=\epsilon(\pi) >0$
such that 

(1) any finite volume K\"ahler manifold 
from $\mathcal M_{n,-4-\epsilon,-1,\pi}$ 
is diffeomorphic to  a scalar multiple of a complex hyperbolic manifold.

(2) any quaternionic-K\"ahler manifold from 
$\mathcal M_{n,-4-\epsilon,-1,\pi}$ 
is diffeomorphic to  a scalar multiple of a quaternionic hyperbolic manifold. 
\end{cor}
\begin{proof}
Arguing by contradiction find a sequence $M_k$ of finite volume K\"ahler
or quaternionic-K\"ahler manifolds with fundamental group isomorphic to
$\pi$ and sectional curvatures within $[-4-1/k,-1]$.
By~\ref{intro: main thm}, we can assume that $(M_k, g_k)$ converges in 
$C^{1,\alpha}$
topology to a $C^{1,\alpha}$-Riemannian manifold $(M,g)$. 
 
First, assume that each $M_k$ is quaternionic-K\"ahler. 
Any quaternionic-K\"ahler is Einstein so $M$ is smooth Einstein manifold
and convergence is in $C^\infty$ topology. Hence, sectional curvatures of
$M_k$ converge to sectional curvatures of
$M$, so $M$ is quarter-pinched, and we are done by~\cite{Yeu}.

Second, assume that each $(M_k, g_k)$ is K\"ahler. 
Adapting the argument in~\cite{Ber} to negative curvature, we deduce that the
holomorphic sectional curvatures of $M_k$ converge to $-1$
uniformly on $M_k$. Look at the "curvature $4$-tensor" 
$R_k^0$ for K\"ahler metric
of holomorphic sectional curvature $-1$ defined in terms of
$g_k$ and almost complex structure $J_k$ of 
$M_k$~\cite[IX.7, right before 7.2]{KN2}.
Since sectional curvature can be written 
in terms of holomorphic sectional curvature~\cite{BG} 
and the curvature $4$-tensor can be written in terms of
sectional curvature, 
the $4$-tensor $R_k$ of $g_k$ is getting close to 
$R_k^0$ uniformly on $M$ when $k\to\infty$.
Taking traces we conclude that Ricci tensor of $g_k$ 
is getting close to $-(n+1)g_k/2$ (see~\cite[IX.7.5]{KN2}).
Now $g_k$ subconverges to a $C^{1,\alpha}$-Riemannian metric on 
a finite volume manifold $M$ hence $-(n+1)g_k/2-Ric(g_k)$
converges to zero. Then the proof of~\cite[theorem 0.4]{Gao}
implies that the limiting metric $g$ is a weak
solution of the Einstein equation, hence $g$ is
a $C^\infty$ Einstein metric.
Also $(M,g)$ has Alexandrov curvature within $[-4,-1]$
hence $sec(M,g)\in [-4,-1]$ and we are done by~\cite{Yeu}. 
\end{proof}

\begin{rmk}
It is not clear whether the assumptions of~\ref{kahler pinching}
are necessary. However, there do exist compact negatively curved
K\"ahler $4$-manifolds which are not homotopy equivalent to
locally symmetric manifolds~\cite{MS}. Also in higher dimensions
there are examples of compact almost quarter pinched 
Riemannian manifolds which are homeomorphic but not diffeomorphic to
complex hyperbolic manifolds~\cite{FJ2}.
\end{rmk}

\small
\bibliographystyle{amsalpha}
\bibliography{xxxfv}

\end{document}